\documentclass[11pt,a4paper]{article}
\usepackage{latexsym}

\title{Flag vectors}

\author{Jonathan Fine\relax
\thanks{203 Coldhams Lane, Cambridge, CB1 3HY, England.  
\quad E-mail: \texttt{j.fine@pmms.cam.ac.uk}}
}
\date{1 October 1998}

\newtheorem{theorem}{Theorem}
\newtheorem{definition}[theorem]{Definition}

\newtheorem{problem}[theorem]{Problem}

\newcommand\bfR{{\bf R}}
\newcommand\fbar{\overline{f}}  
\newcommand\Fbar{\overline{F}}  
\newcommand\ftilde{\widetilde{f}}  
\newcommand\ph[1]{{\bf #1}}

\newcommand\bibrule{\rule{2pc}{0.4pt}}

\begin{document}
\maketitle
\begin{abstract}
\noindent
This paper defines for each object $X$ that can be constructed out of a
finite number of vertices and cells a vector $fX$ lying in a finite
dimensional vector space.  This is the flag vector of $X$.  It is hoped
that the quantum topological invariants of a manifold $M$ can be expressed
as linear functions of the flag vector of the $i$-graph that arises from
any suitable triangulation $T$ of $M$.  Flag vectors are also defined for
finite groups and more generally for $n$-ary relations.  Some problems,
and suggested connections with other constructions, particularly that of
the associahedron and so on, conclude the presentation.
\end{abstract}

\section{Introduction}

This paper defines for each object $X$ that can be constructed out of
vertices and cells a vector $fX$ lying in a finite dimensional vector
space.  This is the flag vector of $X$.  An example that indicates the
importance of this problem follows.

Suppose that $M$ is a compact topological manifold, say of dimension $n$. 
Now let $T$ be a triangulation of $M$.  Because $T$ determines $M$, any
topological invariant of $M$ can be calculated via $T$, at least in
principle.  Considered abstractly, $T$ can be described as follows. First,
let $V$ be the set of points of $M$ that are vertices of cells (in fact
simplices) in the triangulation.  Each $n$-dimensional `triangle' or
simplex in $T$ determines an $(n+1)$-element subset of $V$.  Let $C$ be
the collection of all such subsets of $V$.  Assume also that each cell $c$
in $C$ corresponds to just one simplex of $T$.  It now follows that the
combinatorial object $G=(V,C)$, which is an example of what is known as an
$(n+1)$-graph, determines $T$ and thus $M$ up to equivalence.  This paper
will define a flag vector $fG$ for all $(n+1)$-graphs $G$.

Now suppose that $v=v(M)$ is a numeric (or vector) valued topological
invariant of $n$-manifolds.  For each $(n+1)$-graph $G$ a topological
realisation $|G|$ can be produced, and if $|G|$ is a manifold $M$ then we
will define $v(G)$ to be $v(M)$.  The next three definitions describe the
relationship we seek between the topological invariant $v$ and the flag
vector $f$.

\begin{definition}
Suppose that $v(G_1)=v(G_2)$ whenever $fG_1=fG_2$, for any $G_i$ for which
$|G_i|$ is a manifold.  In that case we will say that $v(G)$ \emph{is a
function of $fG$}.
\end{definition}

\begin{definition}
Suppose that $v$ is a function of $fG$ and in addition that whenever a
linear relation such as
\[
    \lambda_1 fG_1 + \dots + \lambda_r fG_r = 0 \qquad \lambda_i \in \bfR
\]
holds between flag vectors of graphs for which $|G_i|$ is a manifold
$M_i$ then the corresponding relation
\[
    \lambda_1 v(M_1) + \dots + \lambda_r v(M_r) = 0 
\]
holds between the values of $v$, then we will say that \emph{$v(G)$ is a
linear function of $fG$}.
\end{definition}

\begin{definition}
Suppose that $fG$ is given as a point in a vector space $F$, and that $v$
is in the above sense a linear function of $f$.  Any linear function on
$F$ that agrees with $v$ when $|G|$ is a manifold will be called a
\emph{formula for $v$ in terms of $f$}.  (Such a formula will not be
unique, unless $F$ is spanned by the $fG$ for which $|G|$ is a manifold.)
\end{definition}

The methods of quantum topology provide a interesting and steadily growing
collection of topological invariants.  This is an extremely active area of
research.  These invariants are usually discovered family by family, and
in more or less explicit form.  The difficult task is usually to
demonstrate topological invariance.  The flag vector approach to
topological invariants is wholesale rather than retail, or top-down rather
than bottom-up.  It is in part inspired by the Vassiliev theory of knot
invariants, which is similarly a wholesale approach.

Each definition of a flag vector defines a family of topological
invariants, namely those that are linear functions of the flag vector. 
Determining such functions is of course not likely to be so easy. Loosely
speaking, the flag vector is a haystack in which one hopes to find
interesting needles.  For such to be useful, it should be neither too
large nor too small.  One wishes to narrow down the search, without
discarding any needles.

The definitions to be given in this paper will apply not only to
$(n+1)$-graphs but more generally to any object that is built up out of a
finite number of vertices and cells, or can be so represented.  Thus, a
flag vector will be defined for finite groups, where the vertices are the
elements and the cells are the equations $ab=c$ that hold between the
elements.

This paper is organised as follows.  The next two sections define first
the shelling vector and then the flag vector of an $i$-graph, and the
following section describes by means of examples the changes that must be
made to accomodate more general vertex-and-cell objects.  Finally, there
is a summary and the statement of some open problems.

The preprints \cite{bib.JF.QTHGFV,bib.JF.SFV} are perhaps best thought of
as preliminary forms of this paper.  The first does, however, contain
additional material.  This paper can usefully be read in conjunction with
\cite{bib.JF.GFP}, which deals with ordinary or $2$-graphs.

\section{The shelling vector}

Throughout this section and the next $G$ will be an $i$-graph, or in other
words a possibly empty collection $C$ of $i$-element subsets of a vertex
set $V$.  Traditionally, such are called edges.

\begin{definition}
A \emph{shelling} $\sigma$ of $G$ is simply an ordering $v_1$, $\ldots$,
$v_r$ of the vertex set $V$.  The \emph{$j$-th link $L_j$} of the shelling
$\sigma$ is the $(i-1)$-graph consisting of all $(i-1)$-element subsets
$c$ of $v_{j+1}$, $\ldots$, $v_r$ such that $c \cup \{v_j\}$ is a cell of
$G$.
\end{definition}

The definition of the shelling vector is recursive.  Each shelling
$\sigma$ has links $L_i$ which, by assumption, will already have a
shelling vector.  We `multiply' these together, and take the sum over all
shellings.

\begin{definition}
Let $G$ be an $i$-graph, with $i>0$.  The \emph{shelling vector $\ftilde
G$} is the sum
\[
    \sum \nolimits _\sigma \> 
      \ftilde L_1 \otimes \dots \otimes \ftilde L_r
\]
over all shellings of the tensor product of the shelling vectors of the
links.
\end{definition}

\begin{definition}
There is only one set with zero elements, namely the empty set.  Thus, on
any vertex set $V$ there are only two $0$-graphs.  One has no cells, and
the other has the empty set as its only cell.  The shelling vectors of
these graphs will be defined to be the symbols $a$ and $b$ respectively.
\end{definition}

We have now defined a shelling vector for $i$-graphs.  It can be thought
of as a formal sum of words in $a$ and $b$.  (The length of these words
will depend on both the size $r$ of the vertex set and the size $i$ of the
subsets being used.  For example, for $1$-graphs the length is $r$, while
for $2$-graphs the length is $r(r+1)/2$.  In general, the length for
$i$-graphs is the sum of the lengths for the $(i-1)$-graphs that are the
links.)

Almost certainly the shelling vector is too large for our purposes.  If
the following is true then any topological invariant whatsoever of $M$ will
be a linear function of the shelling vector.

\begin{problem}
Suppose $G_1$, $\ldots$, $G_r$ is any collection of $i$-graphs, no two of
which are equivalent.  Are the shelling vectors $\ftilde G_1$, $\ldots$,
$\ftilde G_r$ linearly independent?
\end{problem}

\section{The flag vector}

The flag vector will again be a recursive sum over all shellings, but this
time of products of link contributions.  In the shelling vector, each link
contributed its own shelling vector.  For the flag vector, the
contribution made by a link $L_i$ will be not the whole of the flag vector
$fL_i$ of the link, but some reduction $\fbar L_i$ of this quantity.

\begin{definition}
Let $G$ be an $i$-graph, with $i>0$.  The \emph{flag vector $\ftilde G$}
is the sum
\[
    \sum \nolimits _\sigma \> 
      \fbar L_1 \otimes \dots \otimes \fbar L_r
\]
over all shellings, where the \emph{link contribution $\fbar{L_i}$} will
be defined later.  The flag vector of a $0$-graph is either $a$ or $b$, as
with the shelling vector.
\end{definition}

We will first give the definition of the link contribution, and then we
will motivate it.  Let $F=F_i=F_{i,r}$ be the vector space in which $fG$
naturally lies.  The link contribution $\fbar{L_j}$ will be the residue
$\fbar{L_j}$ of $fL_j$ in a certain quotient $\Fbar$ of $F$.

\begin{definition}
Suppose $G$ is an $i$-graph.  Two distinct cells $c_1$ and $c_2$ of $G$
are \emph{disjoint} if they do not have a vertex in common.  Suppose that
two such cells have been chosen.  Let $G_{+-}$ and $G_{-+}$ denote the
result of removing $c_2$ and $c_1$ respectively from $G$.  Let $G_{--}$
denote the result of removing both $c_1$ and $c_2$, and set $G_{++}$ to be
$G$ itself.  Recall that $fG$ will lie in a vector space $F$.  Define the
\emph{link space} $\Fbar$ to be the quotient of $F$ by the subspace spanned
by
\[
    fG_{++} - fG_{+-} - fG_{-+} + fG_{--}
\]
for all $G$, and all pairs $(c_1,c_2)$ of disjoint cells in $G$.  Now
define the \emph{link contribution} $\fbar G$ of $G$ to be the residue of
$fG$ in $\Fbar$.
\end{definition}

For $0$-graphs the flag vector is by definition the same as the shelling
vector, namely either $a$ or $b$.  For $1$-graphs the links are
$0$-graphs, and as such graphs do not have two cells to be disjoint, the
link contributions are again either $a$ or $b$.  There are up to
equivalence $(r+1)$ distinct $1$-graphs on $r$ vertices, and it is easily
seen that their flag vectors are linearly independent.

For ordinary or $2$-graphs the disjoint cells rule comes in, to reduce the
link contribution from $fL=\ftilde L$ to $\fbar L$.  This is how it goes. 
Just for this paragraph, let $[n]$ denote the link contribution due to an
$n$-celled $1$-graph on some fixed number $r$ of vertices.  The equation
\[
    [n+2] - [n+1] - [n+1] + [n] = 0
\]
follows from the disjoint cell rule, and all relations arise in this way. 
If we write the equation as
\[
    [n+2] - [n+1] = [n+1] - [n]
\]
then its meaning becomes clearer.  The change made by adding a cell to the
link does not depend on the number of cells in the link, at least in the
present situation.  It is easily seen that the vectors
\[
    a = [0] \> , \qquad b = [1] - [0]
\]
provide a basis for the link space, and that $[n]$ is equal to $a+nb$. 
(The symbols $a$ and $b$ are not the same as those used in shelling
vectors.)

We can now motivate the disjoint cell rule for the link contribution.  The
previous paragraph shows that for $2$-graphs the definition in this paper
agrees with that in \cite{bib.JF.GFP}.  The results in \cite{bib.JF.GFP}
indicate that the correct definition has been found, at least for
$2$-graphs, and gives some insight into how the disjoint cell works in
this case.

When the link is a $1$-graph, distinct cells are automatically disjoint. 
Elsewhere \cite{bib.JF.VTRC} the concept of independent regional change
has proved to be useful.  The disjoint cell rule is simply another
application of this principle.

\section{Further examples}

We have defined, for each $i$-graph, a flag vector $fG$.  The same process
can be applied to other combinatorial objects, built up out of vertices
and cells.  For example, to study the topology of oriented manifolds, one
will need to study \emph{oriented} $i$-graphs.  (Such is an $i$-graph,
where each cell has been given an orientation.  An orientation is an
ordering of the vertices, up to an even permutation.)

Once suitable sign conventions have been established, the result of
removing a vertex from an oriented $i$-cell will be an oriented
$(i-1)$-cell.  This works for $i \geq 3$.  For $i=2$ there will be only a
single ordering for the resulting $1$-cell, and so some other convention
must be used instead.  Better however is to change the definition of an
orientation.  Instead, say that an orientation of a cell is a rule that
assigns a sign to each ordering of the vertices, in a manner that respects
even and odd permuations of the vertices.  When this is done, both
$1$-cells and $0$-cells can be oriented.  (In both cases, there is only
one ordering of the vertices available.)

In the unoriented case, the inductive definition of the shelling and flag
vectors was founded on the values $a$ and $b$ for $0$-graphs.  In the
oriented case, there are two possible sorts of $0$-cell, which can
conveniently be denoted by $b_+$ and $b_-$.  These, together with $a$,
will found the inductive definition for oriented $i$-graphs.

Something similar can be done for manifolds with a boundary.  Define an
\emph{$i$-cell with boundary} to be an ordinary $i$-cell, together will a
possibly empty subset of the vertex set, which is the \emph{boundary} of
the cell.  Define an $i$-graph with boundary to be a set of $i$-cells with
boundary, where as ordinary cells the $i$-cells are distinct. Clearly, the
result of removing a vertex from an $i$-cell with boundary will be an
$(i-1)$-cell with boundary.

However, there is more.  We will wish to record whether or not the removed
vertex was a boundary vertex.  (In the oriented case, we used the removed
vertex to choose an orientation for the resulting cell.)  We can record
this as a label attached to the $(i-1)$-cell that results from removing a
vertex from an $i$-cell with boundary.  Successive removal of vertices
will result in a word being written on the label, and so the induction
will be founded on the $0$-cells $a$ and $b_w$, where $a$ is the `empty'
$0$-graph, and where $b_w$ is the `full' $0$-graph whose label $w$ is a
word of length $i$ in say $0$ and $1$, which records the removal of
boundary vertices.

A flag vector for finite groups was promised in the introduction.  The
vertices will be the elements of $G$.  We will think of the group law as
the ternary relations $R=R(x,y,z)$ whose triples $(a,b,c)$ solve the
equation $xy=z$ that represents the group law.  

Loosely speaking, the cells will be the triples $(a,b,c)$ that belong to
$R$.  Given a shelling of the vertices, we will obtain the $i$-th link by
first restricting the relation to the vertices not already removed, and
then taking only those cells that use the $i$-th vertex $v_i$.  However,
because an unordered triple $\{a,b,c\}$ may support several ordered
triples that satisfy $R$, and because a doubleton $\{a,b\}$, or even a
singleton $\{a\}$ may support an ordered triple satisfying $R$, a more
careful approach is needed.

We will instead follow the logic used in \S 2 to define the shelling
vector, but apply it instead to an arbitrary ternary relation $R$.  When a
vertex $v$ is removed from an $i$-graph $G$, the link was defined to be
the $(i-1)$-graph whose cells $c$ become cells of $G$ when $v$ is added to
$c$.  A cell of an $i$-graph is an unordered set of $i$ elements, and
similarly for $(i-1)$ graphs.  Adding a vertex to $c$ is then just a
matter of adding an element to an unordered set.  This is straightforward.

We now apply the same logic to the ternary relation $R$.  Here order and
location are important.  Suppose that $(a,b,c)$ is a triple that satisfies
$R$, with $a$, $b$ and $c$ distinct.  The result of removing $v_i=a$ can
conveniently be denote by $(\ph1,b,c)$, where the $\ph1$ is a
\emph{placeholder} that indicates that a vertex was removed at the first
step.  Similarly, removing $b$ from $(a,b,b)$ will produce
$(a,\ph1,\ph1)$.  The following two definitions are easily seen to be
equivalent.

\begin{definition}
Suppose that $R$ is a ternary relation on a vertex set $V$, and that $a$
is a vertex of $R$.  Now replace $a$ throughout by the placeholder
$\ph1$ in both $R$ and $V$, and discard from $R$ the triples that do not
contain the placeholder $\ph1$.  This new ternary relation $R_a$ is
defined to be the \emph{link} of $R$ at $a$.  (The $i$-th link of a
shelling is much as before $R'_a$, where $R'$ is $R$ restricted to the
vertices that remain before $v_i=a$ is removed.)
\end{definition}

\begin{definition}
The \emph{link} $R_a$ consists of the ordered triples on $V \setminus \{a\}
\cup \{\ph1\}$ that contain at least one placeholder, and which belong to
$R$ when $\ph1$ is replaced by $a$.
\end{definition}

Although the link $R_a$ is still a ternary relation, it is rather
different from its parent $R$ in that its relations all contain the
placeholder $\ph1$, which is not a vertex of $R$.  This must be taken into
account when a second vertex $b$ is removed, to compute the flag vector of
the link.

\begin{definition}
Suppose that $a$ and $b$ are distinct.  The \emph{second-order link}
$R_{ab}$ is obtained from a first order link $R_a$ in the following way. 
In $R_a$ replace $b$ throughout by the placeholder $\ph2$, and discard
from $R_a$ all triples that do not contain at least two placeholders.  The
result is $R_{ab}$.  Similarly, the \emph{third order link} $R_{abc}$ is
obtained from $R_{ab}$ by replacing $c$ by the placeholder $\ph3$ and
discarding from $R_{ab}$ the triples that do not contain at least three
placeholders.
\end{definition}

As with $R_a$ and $R$, so $R_{ab}$ consists of the triples in $V$
$\setminus \{a,b\} \cup \{\ph1,\ph2\}$ with at least two placeholders,
that belong to $R_a$ when $\ph2$ is replaced by $b$.  Much the same holds
for $R_{abc}$ and $R_{ab}$.  (This use of placeholders will define links
not only for ternary relations, but for $n$-ary relations for any $n$, and
for more general objects yet.)

When we get down to the third order link $R_{abc}$, the vertex set $V$
will have disappeared completely.  All that remains will be a ternary
relation $R_{abc}$ on the placeholders $\ph1$, $\ph2$ and $\ph3$.  Not all
such will arise in this way.  For example $(\ph1,\ph3,\ph3)$ is forbidden
because on replacing $\ph3$ by $c$ the result $(\ph1,c,c)$ does not have
at least two placeholders.  The inductive definition of the shelling
vector for ternary relations is now complete.  It is founded on the
ternary relations on $\ph1$, $\ph2$ and $\ph3$ that are not as just
described forbidden.

To go on to define the flag vector we must know what a cell is, so that we
can talk of disjoint pairs of cells.  But first we consider a detail that
arises for all flag vectors, except those of $i$-graphs.

\begin{definition}
Whatever a cell may be, its \emph{support} consists of the set all
vertices of $V$ (but not placeholders) that appear in the cell.
\end{definition}

For $i$-graphs, a cell is determined by its support.  In the oriented
case, and in other situations, there can be several different cells with
the same support set.  Hitherto we have been defining $G_{+-}$ and so on
via the appearance or non-appearance of cells in a disjoint pair.  For
$i$-graphs this was the only choice.

However, for other objects we might wish not to remove cells, but merely
to change them, without altering the support.  It is easily seen that this
additional freedom produces no new relations, in the construction of the
link space $\Fbar$ from the flag vector space $F$ via the disjoint pair of
cells rule.  This is because to change a cell is to first remove it, and
to then replace it by the new value.

\begin{definition}
Let $R'$ be a link for a relation $R$.  For example, $R'$ might be $R_a$
or $R_{ab}$.  A \emph{simple change} consists of the addition to or
removal from $R'$ of an $n$-tuple.  Its \emph{support} is the support of
the cell.
\end{definition}

The next definition is slightly subtle, because for relations the support
of a cell in the link might be empty.  The triples $(\ph1,\ph2,\ph2)$ and
$(\ph1,\ph1,\ph1)$ are examples of this.  Its effect is to group together
such cells in the link into a single compound cell, similar to the
placeholder relations that are used to found the induction.

\begin{definition}
Let $R'$ be as before, and suppose that two simple changes are given,
whose support sets are distinct and disjoint.  As before, this produces
four variants $R'_{++}$, $R'_{+-}$, $R'_{-+}$ and $R'_{--}$ of
$R'=R'_{++}$.  The \emph{disjoint cell rule} defines the \emph{link space}
$\Fbar$ to be the quotient of the vector space $F$ that $fR'$ naturally
lies in by the subspace spanned by
\[
    fR'_{++} - fR'_{+-} - fR'_{-+} + fR'_{--}
\]
for all such $R'$ equipped with such a pair of simple changes.  As before,
the \emph{link contribution} $\fbar {R'}$ of $R'$ is defined to be the
residue of $fR'$ in $\Fbar$.
\end{definition}

This completes the definition of the flag vector for $n$-ary relations,
and group laws $xy=z$ in particular.  It should now be clear how to define
a flag vector for anything that is built out of vertices and cells.

\section{Summary and conclusions}

In this paper we have assigned a flag vector $fG$ to $i$-graphs, groups
and other objects constructed out of vertices and cells, such as $n$-ary
relations.  The example of quantum topology shows the importance of being
able to make such a construction.  In this final section we discuss some
questions whose solution will have some bearing on the fruitfulness of
these definitions.  We begin with quantum topology and the shelling
concept.

\begin{problem}
Can existing quantum topological invariants be expressed as linear
functions of the flag vector?
\end{problem}

Some quantum topology invariants can be computed with the aid of a generic
height function $h$ (a Morse function) on the manifold being studied. 
Such will have a finite number of critical points, each locally equivalent
to a non-degenerate quadratic form.  Now let $T$ be a triangulation of $M$
that is compatible with the height function $h$.  More exactly, adjust $h$
by composing with a monotonic function so that the vertices have as
heights the integers $1$ through to the number $r$ of vertices, and insist
that $h$ is topologically equivalent to its linear interpolation onto the
cells of the the triangulation.

In other words, at least some of the time a shelling $\sigma$ of a
hypergraph $G$ can represent a Morse function $h$ on a manifold $M$.  This
establishes another point of contact between the two theories, and gives
some new insight into the significance of the concept of a shelling. 
Incidentally, in convex polytope theory a shelling of a simplicial
polytope is equivalent (under polarization) to a generic height function
on a simple polytope which in turn, via the moment map, can induce a Morse
function on a projective toric variety.  (The induced function is Morse if
the toric variety is a nonsingular.)

Assume that $G$ and $\sigma$ represent a Morse function $h$ on a manifold
$M$.  In general terms, certain quantum invariants $v$ of $M$ can be
expressed as very special linear combinations of numbers $\lambda_i$ that
can be computed from $G$ and $\sigma$.  On the other hand, each shelling
$\sigma$ of $G$ makes a contribution $f_\sigma G$ to $fG$.  If the
$\lambda_i$ turn out to be linear functions of $f_\sigma G$, then this is
evidence for $v(M)$ being a linear function of $fG$.  For example, if
every shelling $\sigma$ represents a Morse function on $M$, then the
result would follow.

Next we ask to what degree the flag vector distinguishes inequivalent
objects.  The following successively stronger questions are an example of
what we might wish for.

\begin{problem}
Suppose $fG_1=fG_2$.  Does it follow that $G_1$ and $G_2$ are equivalent?
\end{problem}

\begin{problem}
Let $\Delta$ be the convex hull of $fG$, as $G$ runs over a class of
objects.  Are the vectors $fG$ vertices of $\Delta$, and are they distinct
for distinct $G$?
\end{problem}

\begin{problem}
Is there a natural inner product on the space $F$ in which vectors $fG$
(and $\Delta$) lies, such that the $fG$ are distinct and lie on a sphere?
\end{problem}

If the last problem has a positive solution, it is perhaps the easiest way
to resolve the first two.  Although such a result might be thought
unlikely, there are already examples of combinatorially derived polytopes,
whose vertices lie on a sphere.  The permutahedron, the associahedron and
the permuto-associahedron are examples of this. 
Initially~\cite{bib.MK.PA}, the permuto-associahedron was a
combinatorially labelled cell complex whose realization, after some work,
was found to be a sphere. (This result was useful in algebraic topology.) 
Later, it was found that this cell complex could be realized as the
boundary of a convex polytope whose vertices, for a suitably cunning
construction, would lie on a sphere. This phenomenom is at present
unexplained, and is rather strange.  It involves fiber polytopes
\cite{bib.LB-BS.FP}. Lecture~9 of \cite{bib.GZ.LP} is a good first
reference.  The construction is there presented as the problem of
constructing a polytope with prescribed and rather special combinatorics.

This problem can be explored in several ways.  One is to study the flag
vectors of rather special and simple objects, such as binary relations or
partial orders on a small number of vertices.  (A graph is, of course, a
special type of relation among its vertices.)  If the result is a convex
polytope, whose edges and so forth have combinatorial significance, then
we are encouraged.  Another is to define a flag vector for the
combinatorial objects that are the vertices of the associahedron and so
forth.  This will require some thought, for such objects are not at least
on the face of it $n$-ary relations or the like.  If all is well, this
should give an alternative approach to the presently unexplained
construction of the permuto-associahedron.

This paper consists largely of definitions.  Its purpose is to delimit an
area of study, rather than to obtain results in that area.  The flag
vector at present stands somewhat apart from the rest of mathematics, and
other than \cite{bib.JF.GFP} results are not yet available.  A number of
worked examples, and investigation of some of the simpler problems, seems
to be the next step.

Finally, there is another approach.  An object consisting of vertices and
cells can be shelled by removing the vertices, and thus removing the
cells.  This is the path we have followed.  The other approach, described
for graphs at the end of \cite{bib.JF.GFP}, is to shell the object by
removing the \emph{cells} one at a time.  This approach is at the time of
writing completely unexplored.


\begin{thebibliography}{9}

\bibitem{bib.LB-BS.FP}
L.J. Billera and B. Sturmfels, Fiber polytopes,
{\it Annals of Mathematics}, {\bf 135} (1992), 527--549

\bibitem{bib.JF.QTHGFV}
J.Fine, 
On quantum topology, hypergraphs and flag vectors,
preprint q-alg/9708001 \hfil\break  (August 1997)  

\bibitem{bib.JF.SFV}
\bibrule,
On shelling and flag vectors,
preprint q-alg/9710002 (October 1997)

\bibitem{bib.JF.VTRC}
\bibrule,
Vassiliev theory and regional change,
preprint math.QA/9803004 (March 1998)

\bibitem{bib.JF.GFP}
\bibrule,
Graphs, flags and partitions,
preprint math.CO/9809092 (September 1998)

\bibitem{bib.MK.PA}
M. M. Kapranov,
Permuto-Associahedron, MacLane coherence theorem and asymptotic zones for
the KZ equation, 
{\it J.\ Pure and Applied Algebra}, {\bf 85} (1993), 119--142

\bibitem{bib.GZ.LP}
G.M. Ziegler, {\it Lectures on Polytopes},
Springer Verlag, GTM~152 (1995)

\end{thebibliography}
\end{document}